\documentclass[12pt,reqno]{amsart}
\usepackage{amssymb}
\renewcommand*{\leq}{\leqslant} \renewcommand*{\geq}{\geqslant}

\newcommand*{\vare}{\varepsilon} \newcommand*{\varp}{\varphi}

\newcommand*{\mC}{\mathbb C} \newcommand*{\mN}{\mathbb N}
\newcommand*{\mR}{\mathbb R} \newcommand*{\mT}{\mathbb T}
\newcommand*{\mZ}{\mathbb Z}

\newcommand*{\cA}{\mathcal A} \newcommand*{\cB}{\mathcal B}
\newcommand*{\cC}{\mathcal C} \newcommand*{\cG}{\mathcal G}
\newcommand*{\cL}{\mathcal L} \newcommand*{\cM}{\mathcal M}
\newcommand*{\cN}{\mathcal N} \newcommand*{\cO}{\mathcal O}
\newcommand*{\cP}{\mathcal P} \newcommand*{\cT}{\mathcal T}

\newcommand{\cTp}{\cT^{\,\prime}}

\newcommand*{\fd}{\mathfrak d} \newcommand*{\fp}{\mathfrak p}

\newcommand*{\fA}{\mathfrak A} \newcommand*{\fG}{\mathfrak G}
\newcommand*{\fK}{\mathfrak K} \newcommand*{\fN}{\mathfrak N}
\newcommand*{\fO}{\mathfrak O} \newcommand*{\fS}{\mathfrak S}
\newcommand*{\fX}{\mathfrak X}

\newcommand*{\bx}{\overline{x}} \newcommand*{\by}{\overline{y}}
\newcommand*{\bz}{\overline{z}}

\newcommand*{\bX}{\overline{X}}

\newcommand*{\bfX}{\overline{\fX}}

\newcommand*{\hp}{\widehat{p}}

\newcommand*{\hF}{\widehat{F}} \newcommand*{\hM}{\widehat{M}}

\newcommand*{\tF}{\widetilde{F}}

\newcommand*{\tfX}{\widetilde{\fX}}

\newcommand*{\tbeta}{\widetilde{\beta}}

\newcommand*{\bfzero}{\mathbf 0}

\newcommand*{\bfomega}{\text{\boldmath $\omega^2$}}

\newcommand*{\new}{^{\mathrm{new}}}

\newcommand*{\gl}{\mathfrak{gl}} \newcommand*{\GL}{\mathrm{GL}}

\newcommand*{\const}{\mathrm{const}}

\DeclareMathOperator{\codim}{codim} \DeclareMathOperator{\meas}{meas}

\DeclareMathOperator{\Fix}{Fix}

\newtheorem{lem}{Lemma} \newtheorem{thm}{Theorem}

\theoremstyle{definition}
\newtheorem{dfn}{Definition}

\allowdisplaybreaks

\begin{document}

\title[Quasi-periodic perturbations in the reversible context~2]%
{Herman's approach to quasi-periodic perturbations \\ in the reversible KAM context~2}

\author[M.~B.~Sevryuk]{Mikhail B. Sevryuk}

\address{Talroze Institute of Energy Problems of Chemical Physics, The Russia Academy of Sciences, Leninski\u{\i} prospect~38, Bldg.~2, Moscow 119334, Russia}

\email{2421584@mail.ru}

\dedicatory{To the blessed memory of V.~I.~Arnold, a mathematician par excellence}




\begin{abstract}
We revisit non-autonomous systems depending quasi-periodically in time within the reversible context~2 of KAM theory and obtain Whitney smooth families of invariant tori in such systems via Herman's method. The reversible KAM context~2 refers to the situation where the dimension of the fixed point manifold of the reversing involution is less than half the codimension of the invariant torus in question.
\end{abstract}

\maketitle

\baselineskip=17pt

\section{Introduction}\label{intro}

KAM theory founded by the mathematical geniuses A.~N.~Kolmogorov, V.~I.~Arnold, and J.~Moser in 1954--1967 is the theory of quasi-periodic motions in non-integrable dynamical systems. The main ``informal'' conclusion of KAM theory is that Cantor-like families of invariant tori carrying quasi-periodic motions (conditionally periodic motions with incommensurable frequencies) are a generic phenomenon in dynamical systems. However, the properties of such families depend strongly on the phase space structures the system in question is assumed to preserve and on the way the invariant tori relate to these structures. If we confine ourselves with flows on finite-dimensional manifolds, then the following four branches (or, as one often says, \emph{contexts}) of KAM theory are usually considered.

First, KAM theory for Hamiltonian systems, where one should distinguish isotropic (in particular, Lagrangian) invariant tori, coisotropic invariant tori, and so-called atropic invariant tori (i.e., tori that are neither isotropic nor coisotropic). KAM theory for coisotropic and atropic invariant tori can only be developed in the case where the symplectic form $\bfomega$ is not exact \cite{BHS96LNM,BS10}.

Second, KAM theory for reversible systems. Recall that a dynamical system is said to be \emph{reversible} with respect to a smooth involution $G$ of the phase space (a mapping whose square is the identical transformation) if this system is invariant under the transformation $(\fp,t)\mapsto(G\fp,-t)$ where $\fp$ is a point of the phase space and $t$ is the time. In the reversible KAM theory, one always deals with only those tori that are invariant under both the system itself and the reversing involution.

Third, KAM theory for volume preserving systems, where two cases are to be treated separately: invariant tori of codimension one and invariant tori of codimension greater than one.

Fourth, KAM theory for general (dissipative) systems, where the phase space is equipped with no special structure.

In KAM theory, there are also results for some ``exotic'' classes of dynamical systems, for instance, for weakly reversible systems (where the reversing diffeomorphism of the phase space is not assumed to be an involution) \cite{AS86,S86}, for locally Hamiltonian vector fields $V$ (defined by the condition that the $1$-form $i_V\bfomega$ is closed but not necessarily exact, so that the Hamilton function can be multi-valued) \cite{LP05NK,LP07UMZ,S08N}, for conformally Hamiltonian vector fields $V$ (defined by the identity $d(i_V\bfomega)\equiv\eta\bfomega$ with constant $\eta\neq 0$) \cite{CCL13JDE}, for generalized Hamiltonian (or Poisson--Hamilton) systems defined on Poisson manifolds \cite{LY02ETDS,LZH06JMAA} (see \cite{BS10,S08N} for more references), for presymplectic systems (defined in another way on Poisson manifolds where the role of the symplectic form $\bfomega$ is played by a closed degenerate $2$-form with constant rank) \cite{AL12JDDE}, for $b$-Hamiltonian vector fields on the so-called $b$-Poisson (or log-symplectic) manifolds \cite{KMS16JMPA}, or for equivariant vector fields \cite{W10DCDSS}. Here $i_V\bfomega$ is the interior product, or the contraction, of $\bfomega$ with $V$. For recent reviews of various aspects of KAM theory, the reader is referred to the monograph \cite[\S~6.3]{AKN06}, the survey \cite{BS10}, the memoir \cite{GHL14}, and the monograph \cite[Ch.~4]{HCFLM16}. Note that the task of initiating studies of weakly reversible systems was stated by V.~I.~Arnold in \cite{A84}.

In most cases, invariant tori constructed in KAM theory are reducible. The concept of a reducible invariant torus is of principal importance, and we recall here its precise definition.

\begin{dfn}\label{reducible}
Let $\cT$ be an invariant $n$-torus of some flow on an $(n+\fN)$-dimensional manifold. This torus is said to be \emph{reducible} (or \emph{Floquet}) if in a neighborhood of $\cT$, there exists a coordinate frame $x\in\mT^n$, $\fX\in\cO_{\fN}(0)$ in which the torus $\cT$ itself is given by the equation $\{\fX=0\}$ and the dynamical system takes the \emph{Floquet form}
\begin{equation}
\dot{x}=\omega+O(\fX), \quad \dot{\fX}=\Lambda\fX+O_2(\fX)
\label{Floquet}
\end{equation}
with $x$-independent vector $\omega\in\mR^n$ and matrix $\Lambda\in\gl(\fN,\mR)$. The vector $\omega$ (not determined uniquely) is called the \emph{frequency vector} of the torus $\cT$, while the matrix $\Lambda$ (not determined uniquely) is called the \emph{Floquet matrix} of $\cT$, and its eigenvalues are called the \emph{Floquet exponents} of $\cT$. The coordinates $(x,\fX)$ are called the \emph{Floquet coordinates} for $\cT$.
\end{dfn}

In this definition, $\mT^n=(\mR/2\pi\mZ)^n$ is the standard $n$-torus, $\cO_{\fN}(0)$ denotes an unspecified neighborhood of the origin in $\mR^{\fN}$, and $O_2(\fX)$ means $O\bigl(|\fX|^2\bigr)$. We will use similar notation throughout the paper. A crucial property of Cantor families of (reducible) invariant tori in KAM theory is that those families are in fact \emph{Whitney smooth}. This means that although the Floquet coordinates for the tori within a given family are defined a~priori on a certain Cantor-like set, these coordinates can be continued to smooth (say, $C^\infty$) functions defined in an open domain of the appropriate Euclidean space. Whitney differentiability of Cantor families of invariant KAM tori was discovered by V.~F.~Lazutkin around 1970 for the case of invariant curves of area preserving mappings of an annulus (his first well-known paper on the subject was \cite{L71DAN}; for basic references on Whitney smoothness in KAM theory, see \cite{BHS96LNM,BS10}).

As was first pointed out in \cite{BHS96G,BHS96LNM}, the reversible context of KAM theory splits up into two subcontexts with quite different properties, namely, the reversible context~1 and the reversible context~2. Very roughly speaking, consider the following two situations: systems of the form
\begin{equation}
\dot{x}=F(y,\lambda)+O(z), \quad \dot{y}=O(z), \quad \dot{z}=M(y,\lambda)z+O_2(z)
\label{eq1}
\end{equation}
reversible with respect to the involution
\begin{equation}
G_1:(x,y,z)\mapsto(-x,y,Rz),
\label{rev1}
\end{equation}
and systems of the form
\begin{equation}
\dot{x}=F(\lambda)+O(y,z), \quad \dot{y}=\sigma(\lambda)+O(y,z), \quad \dot{z}=M(\lambda)z+O_2(y,z)
\label{eq2}
\end{equation}
reversible with respect to the involution
\begin{equation}
G_2:(x,y,z)\mapsto(-x,-y,Rz).
\label{rev2}
\end{equation}
In both the situations, $x\in\mT^n$, $y\in\mR^m$, $z\in\cO_{2p}(0)$ are the phase space variables ($y\in\cO_m(0)$ in the case of \eqref{eq2}--\eqref{rev2}), $\lambda\in\mR^\kappa$ is an external parameter ($n$, $m$, $p$, $\kappa$ being non-negative integers and $m>0$ in the case of \eqref{eq2}--\eqref{rev2}), $R\in\GL(2p,\mR)$ is an involutive matrix with eigenvalues $1$ and $-1$ of multiplicity $p$ each, $M$ is a $2p\times 2p$ matrix-valued function, and $MR\equiv-RM$. We suppose that in the equation for $\dot{y}$ in~\eqref{eq1}, the terms linear in $z$ are independent of $x$ and, similarly, that in the equation for $\dot{y}$ in~\eqref{eq2}, the terms linear in $y$ and $z$ are independent of $x$.

For each value of $\lambda$, the system~\eqref{eq1} and involution~\eqref{rev1} admit the family $\{y=\const, \, z=0\}$ of reducible invariant $n$-tori carrying conditionally periodic motions with frequency vectors $F(y,\lambda)$, and one is looking for invariant $n$-tori close to $\{y=\const, \, z=0\}$ in small $G_1$-reversible perturbations of family~\eqref{eq1}. This is the prototype of the reversible context~1. Probably the most natural situation where the systems~\eqref{eq1} and their reversible perturbations arise is the search for invariant tori near equilibria of reversible systems (provided that these equilibria are fixed under the reversing involution), see \cite[\S~4.1.3]{BHS96LNM} and \cite{S94} (and references therein).

On the other hand, for each value of $\lambda$ such that $\sigma(\lambda)=0$ (if $\kappa\geq m$ then $\sigma^{-1}(0)$ is generically a $(\kappa-m)$-dimensional surface in $\mR^\kappa$), the system~\eqref{eq2} and involution~\eqref{rev2} admit the reducible invariant $n$-torus $\{y=0, \, z=0\}$ carrying conditionally periodic motions with frequency vector $F(\lambda)$, and one is looking for invariant $n$-tori close to $\{y=0, \, z=0\}$ in small $G_2$-reversible perturbations of family~\eqref{eq2}. This is the prototype of the reversible context~2. The formal definition is as follows.

\begin{dfn}[\cite{BHS96G,BHS96LNM}]\label{firstsecond}
Let $\cT$ be a torus invariant under a system reversible with respect to an involution $G$, and let also $G(\cT)=\cT$. Suppose that the fixed point manifold $\Fix G$ of the involution $G$ is not empty and all the connected components of $\Fix G$ are of the same dimension, so that $\dim\Fix G$ is well-defined. The situation where the inequalities $\frac{1}{2}\codim\cT \leq \dim\Fix G \leq \codim\cT$ hold ($\codim\cT$ being the phase space codimension of $\cT$) is called \emph{the reversible context~1}. The opposite situation where the inequality $\dim\Fix G < \frac{1}{2}\codim\cT$ holds is called \emph{the reversible context~2}.
\end{dfn}

Note that the inequality $\dim\Fix G \leq \codim\cT$ is always valid provided that $\cT$ carries quasi-periodic motions \cite{BHS96G,BHS96LNM,S12MMJ}. For both the families~\eqref{eq1} and~\eqref{eq2}, the codimension $c$ of the invariant $n$-tori in question is equal to $m+2p$. However, $\dim\Fix G_1=m+p\geq c/2$ and $\dim\Fix G_2=p<c/2$ (recall that $m>0$ for \eqref{eq2}--\eqref{rev2}). The drastic differences between the reversible contexts~1 and~2 were discussed in detail in the paper \cite{S11RCD}.

``KAM theory is not only a collection of specific theorems, but rather a methodology, a collection of ideas of how to approach certain problems in perturbation theory connected with `small divisors'\,'' \cite[p.~707]{P01}. This collection has become truly huge by now \cite{AKN06,BHS96LNM,BHTB90,BS10,GHL14,HCFLM16,S86}, and one of its important ingredients is the very powerful method proposed by M.~R.~Herman in 1990 in his talk at the international conference on dynamical systems in Lyons (for a brief written record in the problem of so-called vertically translated $n$-tori in $\mT^n\times\mR$, see \cite[\S~4.6.2]{Y92A}). This method is specifically designed for KAM problems with weak nondegeneracy conditions and consists in the following.

First of all, one proves the so-called ``source'' (or Broer--Huitema--Takens-like) theorem for the context in question. In this theorem, one considers systems depending on external parameters and admitting, in the product of the phase space and the parameter space, a smooth or analytic family of reducible invariant tori. Within this family, the frequencies and Floquet exponents vary in the ``most nondegenerate'' way. The source theorem states that \emph{all} the unperturbed tori with frequencies and Floquet exponents satisfying a suitable Diophantine condition persist under small perturbations of the systems. The corresponding perturbed tori possess \emph{the same} frequency vectors and Floquet matrices and constitute a Whitney smooth family.

Now, to construct invariant tori in systems with degeneracies, one introduces \emph{additional} external parameters to remove all degeneracies. To the new systems, the source theorem can be applied. Finally, one ``extracts'' the desired statement about invariant tori in the original systems from the conclusion of the source theorem, making use of the implicit function theorem and an appropriate (as a rule, rather simple) number-theoretical lemma concerning Diophantine approximations on submanifolds of Euclidean spaces (or, as one says, Diophantine approximations of dependent quantities). The core of Herman's approach is that all the cumbersome and tedious ``KAM machinery'' (homological equations, rapidly convergent infinite sequences of coordinate transformations, etc.)\ is only required to prove the source theorem and is \emph{not needed any longer} to infer various corollaries for systems with degeneracies.

For the isotropic Hamiltonian context, reversible context~1, volume preserving context (for invariant tori of any positive codimension), and dissipative context, the source theorems were presented in the first part (written by H.~W.~Broer, G.~B.~Huitema, and F.~Takens) of the memoir \cite{BHTB90} and in the paper \cite{BH95JDDE}. Some improvements and generalizations are contained in \cite{BCHV09PD,BHN07JDE,W10DCDSS}. From these source theorems, we have deduced various KAM results for systems with degeneracies and complications in the contexts indicated via Herman's method or similar techniques \cite{BHS96G,BHS96LNM,BS10,S95C,S97,S01,S06N,S07DCDS,S07TMIS}. In particular, we have examined invariant tori in systems with very weak (R\"ussmann-like) nondegeneracy conditions \cite{BHS96G,BHS96LNM,BS10,S95C} (see also references therein), the partial preservation of frequencies and Floquet exponents \cite{S06N,S07TMIS}, the so-called excitation of elliptic normal modes (i.e., of purely imaginary Floquet exponents of the unperturbed tori) \cite{BHS96LNM,S95C,S97,S01}, and invariant tori in non-autonomous perturbations of partially integrable systems \cite{S07DCDS}.

Herman's method has been used in KAM theory in many other situations. For instance, it was applied to coisotropic \cite{KP01UMZ,LP05NK,LP07UMZ} and atropic \cite{LP07UMZ} invariant tori of Hamiltonian \cite{KP01UMZ} and locally Hamiltonian \cite{LP05NK,LP07UMZ} systems. In \cite{KX14AMC,WXZ10DCDSB,WXZ17DCDS,ZC10FPTA}, Herman's approach was employed in the case of systems with weak nondegeneracy conditions formulated in terms of the Brouwer topological degree. The papers \cite{KX14AMC,WXZ10DCDSB,WXZ17DCDS} consider invariant tori in the reversible context~1 while the article \cite{ZC10FPTA} treats Hamiltonian systems depending quasi-periodically on time.

Reversible systems and their invariant tori in the framework of context~1 are often encountered (and sometimes quite unexpectedly) in many and various problems of mathematics and physics \cite{A84,AS86,BH95JDDE,D76TAMS,LR98PD,RQ92PR,S86}, so it is not surprising that the reversible context~1 of KAM theory is nearly as developed by now as the isotropic Hamiltonian context. On the other hand, the reversible context~2 remained completely unexplored until 2011 although this context is not only interesting by itself, but also essential for studies of the destruction of unperturbed invariant tori with resonant frequencies in the reversible context~1 \cite{S11RCD,S12IM,S16RCD}. Problem~9 in the list of open problems in KAM theory in \cite{S08N} is ``Develop the reversible KAM theory in context~2''. Up to now, the only works where the reversible KAM context~2 is dealt with have been our papers \cite{S11RCD,S12MMJ,S12IM,S16RCD}. The articles \cite{S11RCD,S12MMJ,S12IM} were based on Moser's modifying terms theory \cite{M67MA}. In the very recent paper \cite{S16RCD}, we have succeeded in obtaining the source theorem for the reversible context~2. Our main tool in \cite{S16RCD} was the BCHV (Broer--Ciocci--Han{\ss}mann--Vanderbauwhede) theorem \cite{BCHV09PD} which concerns a certain particular case of the reversible context~1 (namely, the case of systems~\eqref{eq1} with \emph{singular} matrices $M$). Note that the inference in \cite{S16RCD} of the source theorem for the reversible context~2 from the BCHV theorem is itself similar to Herman's arguments.

Now when the source theorem for the reversible context~2 has become available, it is possible to obtain for this context, employing Herman's approach, analogues of the theorems proven in \cite{BHS96G,BHS96LNM,BS10,S95C,S97,S01,S06N,S07DCDS,S07TMIS} for the ``conventional'' KAM contexts, i.e., the isotropic Hamiltonian context, reversible context~1, volume preserving context, and dissipative context (cf.\ the plan in Section~5 of \cite{S11RCD}). The first step in carrying out this program was made already in \cite{S16RCD} where we considered an analogue of the R\"ussmann nondegeneracy condition for \eqref{eq2}--\eqref{rev2} in the absence of the ``normal'' variable $z$ (i.e., for $p=0$). The present paper is the next step. Namely, we construct reducible invariant $(n+N)$-tori in non-autonomous $G_2$-reversible perturbations of the systems~\eqref{eq2} under the assumption that the perturbation term is quasi-periodic in time $t$ with $N$ basic rationally independent frequencies $\Omega_1,\ldots,\Omega_N$ ($N$ being a positive integer). If $\dot{\fp}=V(\fp,\lambda)$ is an abbreviated notation for~\eqref{eq2} with $\fp=(x,y,z)$, then the perturbed systems have the form
\[
\dot{\fp}=V(\fp,\lambda)+\cP(\fp,\lambda,\Omega_1t,\ldots,\Omega_Nt),
\]
where the function $\cP=\cP(\fp,\lambda,X_1,\ldots,X_N)$ is $2\pi$-periodic in each of the variables $X_1,\ldots,X_N$. The natural reformulation of the problem is to look for invariant $(n+N)$-tori in the corresponding autonomous systems
\begin{equation}
\dot{\fp}=V(\fp,\lambda)+\cP(\fp,\lambda,X), \quad \dot{X}=\Omega
\label{fpVcP}
\end{equation}
reversible with respect to the involution $\fG:(x,y,z,X)\mapsto(-x,-y,Rz,-X)$ where $X\in\mT^N$. The analogous problems for the ``conventional'' KAM contexts were treated by Herman's method in the paper \cite{S07DCDS} where also an extensive general bibliography on KAM theory for non-autonomous flows is presented.

In fact, reducible invariant $(n+N)$-tori in $\fG$-reversible systems~\eqref{fpVcP} were obtained in \cite{S12IM}. However, in \cite{S12IM}, we only constructed $(\kappa-n-m-\nu)$-parameter analytic families of such tori where $\nu$ is the number of the eigenvalues of the matrix $M(\lambda)$, for any fixed value of $\lambda$, with positive imaginary parts and non-negative real parts (to be more precise, $\nu=\nu_2+\nu_3$ in the notation of Definition~\ref{fS} below). The number $\kappa$ of external parameters $\lambda_1,\ldots,\lambda_\kappa$ was assumed to be no less than $n+m+\nu$. Moreover, the frequencies and positive imaginary parts of the Floquet exponents of the perturbed invariant $(n+N)$-tori in \cite{S12IM} are the same for all the tori in the given family. The result of the present paper (Theorem~\ref{thmain} below) is much stronger. First, we show that it suffices to require $\kappa\geq m+1$. To be more precise, one needs $m$ external parameters $\sigma_1,\ldots,\sigma_m$ to overcome a drift along the variable $y$ (such a drift is automatically precluded by the $G_1$-reversibility but not by the $G_2$-reversibility) and one more external parameter $\mu$ (which can be just one-dimensional) to control the frequencies and Floquet exponents of the tori sought for. Second, we obtain $(\kappa-m)$-parameter Whitney smooth families of invariant $(n+N)$-tori. Note that the techniques of \cite{M67MA} (the paper \cite{S12IM} is based on) do not enable one to construct Whitney smooth families of invariant tori.

Like in \cite{S07DCDS}, we confine ourselves with analytic systems but there is no doubt that the result can be carried over to Gevrey regular, $C^\infty$-smooth, or finitely differentiable systems. Similarly, the families of analytic perturbed invariant tori in Theorems~\ref{thmain} and~\ref{thsource} below are claimed to be $C^\infty$-smooth in the sense of Whitney, but most probably these families are Gevrey regular in the sense of Whitney (cf.\ \cite{W10DCDSS}).

Note that although the overwhelming majority of works on invariant tori in non-autonomous systems within KAM theory are confined with systems depending on time periodically or quasi-periodically, there are also some results in KAM theory that concern systems depending on time in a more or less arbitrary way, see \cite{FW16} and references therein.

The paper is organized as follows. In Section~\ref{prelim}, we explain all the notation and present the Diophantine lemma to be used in Herman's procedure. The main result of the paper is stated in Section~\ref{mainres}. In Section~\ref{source}, we give a precise formulation of the source theorem for the reversible context~2 \cite{S16RCD} in the form we need. A proof of the main result is expounded in Section~\ref{proof}.

The paper is dedicated to the memory of V.~I.~Arnold, one of the greatest scholars of the second half of the 20th century and, in particular, one of the founders of KAM theory. Many branches of modern mathematics would be inconceivable without his fundamental contributions, both in the form of specific theorems and in the form of a powerful new ideology. The contemporary mathematical culture has benefitted tremendously from his ingenious books and his passionate essays concerning science and education. Many mathematicians (including the author) owe him a huge debt of gratitude for his generous guidance and help throughout their careers.

\section{Preliminaries}\label{prelim}

Let $\mN$ be the set of positive integers and let $\mZ_+=\mN\cup\{0\}$. We will denote by $|{\cdot}|$ the $\ell_1$-norm of vectors in $\mC^s$, by $\|{\cdot}\|$ the $\ell_2$-norm of vectors in $\mR^s$, and by $\langle{\cdot},{\cdot}\rangle$ the inner product of two vectors in $\mR^s$. The Lebesgue measure in $\mR^s$ will be denoted by $\meas_s$. We will adopt the standard multi-index notation
\[
q!=q_1!q_2!\cdots q_s!, \quad
\mu^q=\mu_1^{q_1}\mu_2^{q_2}\cdots\mu_s^{q_s}, \quad
D_\mu^qF=\frac{\partial^{|q|}F}{\partial\mu_1^{q_1}\partial\mu_2^{q_2}\cdots\partial\mu_s^{q_s}},
\]
where $q\in\mZ_+^s$, $\mu\in\mR^s$, and $F$ is a (vector-valued) function $C^{|q|}$-smooth in $\mu$ (however, this notation will only be used in Definition~\ref{nondegenerate} below). The expression $\cO_s(\mu^0)$ will denote an unspecified neighborhood of a point $\mu^0\in\mR^s$. If $d\in\mN$ and $x,y,z,\ldots$ are certain variables, we will write $O_d(x,y,z,\ldots)$ instead of $O\bigl(|x|^d+|y|^d+|z|^d+\cdots\bigr)$ where $O$ is the standard ``big-O'' symbol. Of course, we will write $O({\cdot})$ instead of $O_1({\cdot})$. Instead of $\{0\}$ with $0\in\mR^s$, we will sometimes write $\{0\in\mR^s\}$.

For any $\mR^s$-valued function $a$ of angular variables $x\in\mT^n$ (and, maybe, of some other variables), the notation $\{a\}_x$ will denote the averaging of $a$ over $\mT^n\ni x$.

The $m\times m$ zero matrix will be denoted by $\bfzero_m$ and the space of $n\times m$ real matrices by $\mR^{n\times m}$, so that $\gl(n,\mR)=\mR^{n\times n}$. If $\cA$ and $\cB$ are two square matrices then $\cA\oplus\cB$ will denote the block diagonal matrix with blocks $\cA$ and $\cB$. Let $R\in\GL(p+\hp,\mR)$ be an involutive matrix with eigenvalue $1$ of multiplicity $p$ and eigenvalue $-1$ of multiplicity $\hp$. One says that a matrix $M\in\gl(p+\hp,\mR)$ anti-commutes with $R$, or is \emph{infinitesimally reversible} with respect to $R$, if $MR=-RM$. If this is the case then the eigenvalues of $M$ come in pairs $(\lambda,-\lambda)$, and if $\hp\neq p$ then $0$ is an eigenvalue of $M$ of multiplicity at least $|\hp-p|$ \cite{H96JDE,S86,S91TSP,S93CJM}. If $\hp=p$ and the spectrum of an infinitesimally $R$-reversible matrix $M$ is simple (a generic case), then $M$ is non-singular (because otherwise $0$ would be an eigenvalue of $M$ of multiplicity at least $2$).

\begin{dfn}\label{fS}
Let a matrix $M\in\GL(2p,\mR)$ anti-commute with an involutive $2p\times 2p$ matrix with eigenvalues $1$ and $-1$ of multiplicity $p$ each. We write that the spectrum of $M$ has the form $\fS(\nu_1,\nu_2,\nu_3;\alpha,\beta)$ where $\nu_1,\nu_2,\nu_3\in\mZ_+$, $\nu_1+\nu_2+2\nu_3=p$, and $\alpha\in\mR^{\nu_1+\nu_3}$, $\beta\in\mR^{\nu_2+\nu_3}$ are two vectors with positive components, if the eigenvalues of $M$ have the form
\begin{gather*}
\pm\alpha_1,\ldots,\pm\alpha_{\nu_1}, \qquad \pm i\beta_1,\ldots,\pm i\beta_{\nu_2}, \\
\pm\alpha_{\nu_1+1}\pm i\beta_{\nu_2+1},\ldots,\pm\alpha_{\nu_1+\nu_3}\pm i\beta_{\nu_2+\nu_3}.
\end{gather*}
\end{dfn}

Recall that a $C^1$-smooth mapping $F:\cM^m\to\cN^n$ of an $m$-dimensional manifold $\cM$ to an $n$-dimensional manifold $\cN$ is said to be \emph{submersive} at a point $\mu\in\cM$, if $m\geq n$ and the rank of the differential of $F$ is equal to $n$ at $\mu$. If this is the case then $F$ is also submersive at any point $\mu'\in\cM$ sufficiently close to $\mu$.

\begin{dfn}[\cite{S06N,S07DCDS,S07TMIS}]\label{Diophantine}
Let $n,\nu\in\mZ_+$. Given $\tau\geq 0$, $\gamma>0$, and $L\in\mN$, a pair of vectors $F\in\mR^n$, $\beta\in\mR^\nu$ is said to be \emph{affinely $(\tau,\gamma,L)$-Diophantine}, if the inequality
\[
\bigl|\langle F,k\rangle + \langle\beta,l\rangle\bigr| \geq \gamma|k|^{-\tau}
\]
holds for any $k\in\mZ^n\setminus\{0\}$ and $l\in\mZ^\nu$ such that $|l|\leq L$.
\end{dfn}

Clearly, if $n\in\mN$ and a pair of vectors $F\in\mR^n$, $\beta\in\mR^\nu$ is affinely $(\tau,\gamma,L)$-Diophantine, then the vector $F$ is $(\tau,\gamma)$-Diophantine in the usual sense, so that $\tau\geq n-1$.

\begin{dfn}[\cite{S06N,S07DCDS,S07TMIS}]\label{nondegenerate}
Let $s\in\mN$ and $n,\nu\in\mZ_+$. Let $\fK\subset\mR^s$ be an open domain and let $Q,L\in\mN$. Consider a pair of $C^Q$-smooth mappings $F:\fK\to\mR^n$, $\beta:\fK\to\mR^\nu$. If $n>0$, introduce the notation
\[
\rho^Q(\mu) = \min_{\|e\|=1} \max_{J=1}^Q J! \max_{\|u\|=1}
\left|\sum_{|q|=J} \bigl\langle D_\mu^qF(\mu),e\bigr\rangle \frac{u^q}{q!}\right|
\]
($q\in\mZ_+^s$, $e\in\mR^n$, $u\in\mR^s$) for $\mu\in\fK$. If $\nu>0$, introduce the notation
\[
\Xi_l^Q(\mu) = \max_{J=1}^Q J! \max_{\|u\|=1}
\left|\sum_{|q|=J} \bigl\langle D_\mu^q\beta(\mu),l\bigr\rangle \frac{u^q}{q!}\right|
\]
($q\in\mZ_+^s$, $u\in\mR^s$) for $\mu\in\fK$, $l\in\mZ^\nu$. The pair of mappings $F$, $\beta$ is said to be \emph{affinely $(Q,L)$-nondegenerate} at a point $\mu\in\fK$ if one of the following four conditions is satisfied.

1) $n>0$, $\nu>0$, $\rho^Q(\mu)>0$, and
\[
\max_{1\leq |q|\leq Q} \Bigl| \bigl\langle D_\mu^qF(\mu),k\bigr\rangle + \bigl\langle D_\mu^q\beta(\mu),l\bigr\rangle \Bigr| > 0
\]
($q\in\mZ_+^s$) for all $k\in\mZ^n$ and $l\in\mZ^\nu$ such that $1\leq |l|\leq L$ and $\|k\| \leq \Xi_l^Q(\mu)\big/\rho^Q(\mu)$.

2) $n>0$, $\nu=0$, and $\rho^Q(\mu)>0$.

3) $n=0$, $\nu>0$, and $\Xi_l^Q(\mu)>0$ for all $l\in\mZ^\nu$ such that $1\leq |l|\leq L$.

4) $n=\nu=0$.
\end{dfn}

Note that for any (vector-valued) $C^J$-smooth function $H$ defined in $\fK\subset\mR^s$ ($J\in\mZ_+$) and any $\mu\in\fK$, $u\in\mR^s$ one has
\[
J! \sum_{|q|=J} D_\mu^qH(\mu)\frac{u^q}{q!} =
\left. \frac{d^J}{dt^J}H(\mu+tu) \right|_{t=0}
\]
($q\in\mZ_+^s$). The inequality $\rho^Q(\mu)>0$ (for $n>0$) means that the collection of all the $\binom{s+Q}{s}-1$ partial derivatives of $F$ at $\mu$ of all the orders from $1$ to $Q$ spans $\mR^n$, i.e., the linear hull of these derivatives is $\mR^n$ (a R\"ussmann-type property). The inequality $\Xi_l^Q(\mu)>0$ (for $\nu>0$ and some $l\in\mZ^\nu\setminus\{0\}$) means that at least one of the $\binom{s+Q}{s}-1$ partial derivatives of $\beta$ at $\mu$ of all the orders from $1$ to $Q$ is not orthogonal to $l$. Obviously, if a pair of mappings $F$, $\beta$ is affinely $(Q,L)$-nondegenerate at a point $\mu\in\fK$, then it is affinely $(Q,L)$-nondegenerate at any point $\mu'\in\fK$ sufficiently close to $\mu$. Note however that in this statement, it is essential that the inequality $\|k\| \leq \Xi_l^Q(\mu)\big/\rho^Q(\mu)$ in the first condition of Definition~\ref{nondegenerate} is not strict.

\begin{lem}[\cite{S07DCDS}]\label{lmDioph}
Let $s\in\mN$ and $n,\nu\in\mZ_+$. Let $\fK\subset\mR^s$ be an open domain and $K\subset\fK$ a subset of $\fK$ diffeomorphic to a closed $s$-dimensional ball. Let also $Q,L\in\mN$. Finally, let a pair of $C^Q$-smooth mappings $F:\fK\to\mR^n$, $\beta:\fK\to\mR^\nu$ be affinely $(Q,L)$-nondegenerate at each point of $K$. Then

\textup{1)} there exists a number $\delta>0$ and

\textup{2)} for every $N\in\mN$, $\tau_\ast\geq N-1$, $\gamma_\ast>0$, $\vare\in(0,1)$ and every $\tau$ such that $\tau>(n+N)Q$ and $\tau\geq\tau_\ast$, there exists a number $\gamma=\gamma(\vare,\tau,\gamma_\ast)>0$ such that the following holds. Let $\tF:\fK\to\mR^n$ and $\tbeta:\fK\to\mR^\nu$ be any $C^Q$-smooth mappings such that all the partial derivatives of each component of $\tF-F$ and $\tbeta-\beta$ of any order from $1$ to $Q$ are smaller than $\delta$ in absolute value everywhere in $\fK$. Let $\Omega\in\mR^N$ be an arbitrary $(\tau_\ast,\gamma_\ast)$-Diophantine vector. Then the Lebesgue measure of the set of those points $\mu\in K$ for which the pair of vectors
\begin{equation}
\bigl(\tF(\mu),\Omega\bigr)\in\mR^{n+N}, \quad \tbeta(\mu)\in\mR^\nu
\label{vectors1}
\end{equation}
is affinely $(\tau,\gamma,L)$-Diophantine, is greater than $(1-\vare)\meas_sK$.
\end{lem}

This Diophantine lemma (Lemma~5.1 in \cite{S07DCDS}) is a particular case of Lemma~3.3 in \cite{S06N} which in turn is a particular case of Lemma~2.3 in \cite{S07TMIS}. The importance of studying Diophantine approximations of dependent quantities for many problems in mathematics and mathematical physics involving small divisors was first emphasized by V.~I.~Arnold in 1968 in his talk ``Problems of Diophantine approximations in analysis'' at the symposium on number theory in Vladimir, Russia.

\section{The Main Result}\label{mainres}

Before presenting a rigorous formulation of our main result, we would like to discuss in a more relaxed manner what $G_2$-reversible systems~\eqref{eq2} could look like. The perturbation term $\cP(\fp,\lambda,X)$ in~\eqref{fpVcP} is of course assumed to be an \emph{arbitrary} function subject to the appropriate reversibility and smallness conditions (recall that $\fp=(x,y,z)$). On the other hand, we will suppose that the unperturbed systems $\dot{\fp}=V(\fp,\lambda)$ are ``integrable'' in the sense that they are $\mT^n$-equivariant, i.e., $V$ is independent of $x$. This technical restriction allows one to apply the source theorem \cite{S16RCD} for the reversible context~2. Thus, the unperturbed systems~\eqref{eq2} in question have the form
\begin{equation}
\begin{aligned}
\dot{x} &= F(\lambda)+O(y,z), \\
\dot{y} &= \sigma(\lambda)+Y(\lambda)y+Z(\lambda)z+O_2(y,z), \\
\dot{z} &= M(\lambda)z+O_2(y,z),
\end{aligned}
\label{add1}
\end{equation}
where the $O(\cdot)$ terms are $x$-independent while $Y$ and $Z$ are $m\times m$ and $m\times 2p$ matrix-valued functions, respectively. The $G_2$-reversibility of~\eqref{add1} implies that $Y\equiv 0$ and $ZR\equiv Z$. We will assume that the matrix $M(\lambda)$ governing the linear behavior of the ``normal'' variable $z$ is non-singular for each $\lambda$. Then for $y'=y-ZM^{-1}z$ one has $\dot{y}'=\sigma(\lambda)+O_2(y',z)$. Moreover, in the new coordinate frame $(x,y',z)$, the involution~\eqref{rev2} retains its form $G_2:(x,y',z)\mapsto(-x,-y',Rz)$ because
\[
-y'=-y+ZM^{-1}z=-y-ZM^{-1}Rz
\]
(the identities $MR\equiv-RM$ and $ZR\equiv Z$ imply that $ZM^{-1}R\equiv-ZM^{-1}$). Thus, one may suppose without loss of generality that the term $Zz$ is absent in~\eqref{add1}.

The parameter $\lambda\in\mR^\kappa$ in~\eqref{add1} is to be considered near the set $\sigma^{-1}(0)\subset\mR^\kappa$. We will assume that $s=\kappa-m\geq 1$ and the mapping $\lambda\mapsto\sigma(\lambda)$ is submersive at a certain point $\lambda_\star\in\sigma^{-1}(0)$ (this is an essential hypothesis), so that $\sigma^{-1}(0)$ is a smooth $s$-dimensional surface near $\lambda_\star$. One can choose a new coordinate frame in the parameter space $\mR^\kappa$ near $\lambda_\star$ in the form $(\sigma,\mu)$ where $\sigma\in\cO_m(0)$ and $\mu\in\cO_s(0)$, the point $(\sigma=0, \, \mu=0)$ corresponding to $\lambda=\lambda_\star$. Then the systems~\eqref{add1} take the form
\begin{align*}
\dot{x} &= \hF(\sigma,\mu)+\xi(y,z,\sigma,\mu), \\
\dot{y} &= \sigma+\eta(y,z,\sigma,\mu), \\
\dot{z} &= \hM(\sigma,\mu)z+\zeta(y,z,\sigma,\mu),
\end{align*}
where $\xi=O(y,z)$, $\eta=O_2(y,z)$, $\zeta=O_2(y,z)$. We will write
\begin{equation}
\hF(\sigma,\mu)=F(\mu)+\Delta(\sigma,\mu), \quad \hM(\sigma,\mu)=M(\mu)+\Pi(\sigma,\mu),
\label{add2}
\end{equation}
where $\Delta=O(\sigma)$, $\Pi=O(\sigma)$, $MR\equiv-RM$, and $\Pi R\equiv-R\Pi$. Of course, one should not confuse the functions $F(\mu)=\hF(0,\mu)$, $M(\mu)=\hM(0,\mu)$ in~\eqref{add2} and the functions $F$, $M$ in~\eqref{add1}. Following \cite{S16RCD}, we will set $\zeta=O_2(y,z,\sigma)$ instead of $\zeta=O_2(y,z)$ and ``incorporate'' the term $\Pi(\sigma,\mu)z$ into $\zeta$, cf.\ equations~(2.2) in \cite{S16RCD}.

Thus, the final form of the integrable unperturbed systems in our theorem is
\begin{equation}
\begin{aligned}
\dot{x} &= F(\mu)+\Delta(\sigma,\mu)+\xi(y,z,\sigma,\mu), \\
\dot{y} &= \sigma+\eta(y,z,\sigma,\mu), \\
\dot{z} &= M(\mu)z+\zeta(y,z,\sigma,\mu),
\end{aligned}
\label{add3}
\end{equation}
where $\Delta=O(\sigma)$, $\xi=O(y,z)$, $\eta=O_2(y,z)$, $\zeta=O_2(y,z,\sigma)$. For $\sigma=0$ and any $\mu$, the system~\eqref{add3} admits the reducible invariant $n$-torus $\{y=0, \, z=0\}$ with frequency vector $F(\mu)\in\mR^n$ and Floquet matrix $\bfzero_m\oplus M(\mu)\in\gl(m+2p,\mR)$. This torus is also invariant under the reversing involution~\eqref{rev2}.

Another point we would like to recall is how the perturbed Floquet coordinates in KAM theory can be most simply expressed in terms of the unperturbed ones (this aspect of the theory is standard but rarely expounded in detail). Consider the general setup of Definition~\ref{reducible}. Suppose that a small perturbation of the system~\eqref{Floquet} has a reducible invariant $n$-torus $\cTp$ close to the unperturbed torus $\cT=\{\fX=0\}$. Assume also that the Floquet coordinates for $\cTp$ can be chosen to be close to the Floquet coordinates $(x,\fX)$ for $\cT$. The torus $\cTp$ is given in the parametric form by the equations $x=\bx+a(\bx)$, $\fX=\fA^0(\bx)$ where $\bx\in\mT^n$, $\dot{\bx}=\omega'$, the functions $a:\mT^n\to\mR^n$ and $\fA^0:\mT^n\to\mR^{\fN}$ are small, and $\omega'\in\mR^n$ is the frequency vector of $\cTp$ ($\omega'$ being close to $\omega$). The relations $x=\bx+a(\bx)$, $\fX=\tfX+\fA^0(\bx)$ with $\tfX\in\cO_{\fN}(0)$ determine the new coordinate frame $(\bx,\tfX)$ around $\cTp$. In the coordinates $(\bx,\tfX)$, the torus $\cTp$ is given by the equation $\{\tfX=0\}$ while the perturbed equations of motion take the form
\[
\dot{\bx}=\omega'+O(\tfX), \quad \dot{\tfX}=\bigl(\Lambda+\fd(\bx)\bigr)\tfX+O_2(\tfX),
\]
where the function $\fd:\mT^n\to\gl(\fN,\mR)$ is small. Since the invariant torus $\cTp$ is assumed to be reducible, there exists a coordinate transformation $\tfX=\bfX+\fA^1(\bx)\bfX$ with $\bfX\in\cO_{\fN}(0)$ that casts the linear equation $\dot{\tfX}=\bigl(\Lambda+\fd(\bx)\bigr)\tfX$ (for $\dot{\bx}=\omega'$) into an equation $\dot{\bfX}=\Lambda'\bfX$ with a \emph{constant} matrix $\Lambda'\in\gl(\fN,\mR)$ (the Floquet matrix of $\cTp$) close to $\Lambda$. Here $\fA^1:\mT^n\to\gl(\fN,\mR)$ is a small function. The relations
\begin{equation}
x=\bx+a(\bx), \quad \fX=\bfX+\fA^0(\bx)+\fA^1(\bx)\bfX
\label{add4}
\end{equation}
determine the Floquet coordinates $(\bx,\bfX)$ for $\cTp$. In these coordinates, $\cTp=\{\bfX=0\}$ and the perturbed equations of motion take the Floquet form $\dot{\bx}=\omega'+O(\bfX)$, $\dot{\bfX}=\Lambda'\bfX+O_2(\bfX)$, cf.~\eqref{Floquet}. Note that in~\eqref{add4}, $x$ is independent of $\bfX$ whereas $\fX$ depends on $\bfX$ in an \emph{affine} way. The terms $a(\bx)$ and $\fA^0(\bx)$ are responsible for the torus $\cTp$ itself while the term $\fA^1(\bx)\bfX$ is responsible for the variational equation along $\cTp$.

Now we are in a position to state formally our main theorem. Recall also that the meaning of the notation $\fS(\nu_1,\nu_2,\nu_3;\alpha,\beta)$ for matrix spectra is explained in Definition~\ref{fS}.

Let $n\in\mZ_+$, $m\in\mN$, $p\in\mZ_+$, $N\in\mN$, $s\in\mN$. Consider an analytic $(m+s)$-parameter family of analytic differential equations
\begin{equation}
\begin{aligned}
\dot{x} &= F(\mu)+\Delta(\sigma,\mu)+\xi(y,z,\sigma,\mu)+f(x,y,z,\sigma,\mu,X), \\
\dot{y} &= \sigma+\eta(y,z,\sigma,\mu)+g(x,y,z,\sigma,\mu,X), \\
\dot{z} &= M(\mu)z+\zeta(y,z,\sigma,\mu)+h(x,y,z,\sigma,\mu,X), \\
\dot{X} &= \Omega,
\end{aligned}
\label{eqpert}
\end{equation}
where $x\in\mT^n$, $y\in\cO_m(0)$, $z\in\cO_{2p}(0)$, $X\in\mT^N$ are the phase space variables, $\sigma\in\cO_m(0)$ and $\mu\in\cO_s(0)$ are external parameters, $M$ is a $2p\times 2p$ matrix-valued function, $\Delta=O(\sigma)$, $\xi=O(y,z)$, $\eta=O_2(y,z)$, $\zeta=O_2(y,z,\sigma)$, and the vector $\Omega\in\mR^N$ is $(\tau_\ast,\gamma_\ast)$-Diophantine with some constants $\tau_\ast\geq N-1$ and $\gamma_\ast>0$. The functions $F$, $\Delta$, $M$, $\xi$, $\eta$, $\zeta$ are supposed to be fixed whereas the summands $f$, $g$, $h$ are small perturbation terms. Let the systems~\eqref{eqpert} be reversible with respect to the phase space involution
\begin{equation}
\fG:(x,y,z,X)\mapsto(-x,-y,Rz,-X),
\label{invol}
\end{equation}
where $R\in\GL(2p,\mR)$ is an involutive matrix with eigenvalues $1$ and $-1$ of multiplicity $p$ each, $M(\mu)R\equiv-RM(\mu)$, and the spectrum of $M(0)$ is simple. One may assume that the spectrum of $M(\mu)$ is simple for each $\mu$ and has the form $\fS\bigl(\nu_1,\nu_2,\nu_3;\alpha(\mu),\beta(\mu)\bigr)$ where $\nu_1+\nu_2+2\nu_3=p$. Introduce the notation $\nu=\nu_2+\nu_3\in\mZ_+$.

\begin{thm}\label{thmain}
Suppose that the pair of mappings $F:\cO_s(0)\to\mR^n$, $\beta:\cO_s(0)\to\mR^\nu$ is affinely $(Q,2)$-nondegenerate at $0$ for some $Q\in\mN$. Then there exists a closed $s$-dimensional ball $\Gamma\subset\mR^s$ centered at the origin and such that the following holds. For every complex neighborhood
\begin{equation}
\cC\subset(\mC/2\pi\mZ)^n\times\mC^{2m+2p+s}\times(\mC/2\pi\mZ)^N
\label{neigh}
\end{equation}
of the set
\begin{equation}
\mT^n\times\{0\in\mR^m\}\times\{0\in\mR^{2p}\}\times\{0\in\mR^m\}\times\Gamma\times\mT^N,
\label{set}
\end{equation}
every $\cL\in\mN$, $\vare_1>0$, $\vare_2\in(0,1)$, and every $\tau$ such that $\tau>(n+N)Q$ and $\tau\geq\tau_\ast$, there are numbers $\delta>0$ and $\gamma\in(0,\gamma_\ast]$ with the following properties.

Suppose that the perturbation terms $f$, $g$, $h$ in~\eqref{eqpert} can be holomorphically continued to the neighborhood $\cC$ and $|f|<\delta$, $|g|<\delta$, $|h|<\delta$ in $\cC$. Then there exist a set $\cG\subset\Gamma$, a function $\Theta:\Gamma\to\mR^m$, and a change of variables
\begin{equation}
\begin{aligned}
x &= \bx+a(\bx,X,\mu), \\
y &= \by+b^0(\bx,X,\mu)+b^1(\bx,X,\mu)\by+b^2(\bx,X,\mu)\bz, \\
z &= \bz+c^0(\bx,X,\mu)+c^1(\bx,X,\mu)\by+c^2(\bx,X,\mu)\bz, \\
X &= X
\end{aligned}
\label{change}
\end{equation}
\textup{(}cf.~\eqref{add4}\textup{)} for each $\mu\in\Gamma$ with $\bx\in\mT^n$, $\by\in\cO_m(0)$, $\bz\in\cO_{2p}(0)$ such that the following holds.

\textup{i)} The coefficients $a$, $b^0$, $b^1$, $b^2$, $c^0$, $c^1$, $c^2$ in~\eqref{change} are mappings ranging in $\mR^n$, $\mR^m$, $\gl(m,\mR)$, $\mR^{m\times 2p}$, $\mR^{2p}$, $\mR^{2p\times m}$, $\gl(2p,\mR)$, respectively. These mappings are analytic in $(\bx,X)$ and $C^\infty$-smooth in $\mu$. All the partial derivatives of each component of these mappings of any order from $0$ to $\cL$ are smaller than $\vare_1$ in absolute value everywhere in $\mT^{n+N}\times\Gamma$. The function $\Theta$ is $C^\infty$-smooth, and all the partial derivatives of each component of $\Theta$ of any order from $0$ to $\cL$ are smaller than $\vare_1$ in absolute value everywhere in $\Gamma$. Moreover, $\{a\}_{(\bx,X)}\equiv 0$.

\textup{ii)} For each $\mu\in\Gamma$, the change of variables~\eqref{change} commutes with the involution~\eqref{invol} in the sense that in the new variables $(\bx,\by,\bz,X)$, the involution $\fG$ takes the form
\[
\fG:(\bx,\by,\bz,X)\mapsto(-\bx,-\by,R\bz,-X).
\]

\textup{iii)} $\meas_s\cG > (1-\vare_2)\meas_s\Gamma$.

\textup{iv)} For any point $\mu\in\cG$, the system~\eqref{eqpert} with $\sigma=\Theta(\mu)$ after the coordinate transformation~\eqref{change} takes the form
\begin{equation}
\dot{\bx}=\omega'+O(\by,\bz), \quad \dot{\by}=O_2(\by,\bz), \quad \dot{\bz}=M'\bz+O_2(\by,\bz), \quad \dot{X}=\Omega
\label{goal}
\end{equation}
with $\omega'\in\mR^n$ and $M'\in\GL(2p,\mR)$.

\textup{v)} The $2p\times 2p$ matrix $M'$ satisfies the condition $M'R=-RM'$. The spectrum of $M'$ is simple and has the form $\fS(\nu_1,\nu_2,\nu_3;\alpha',\beta')$. The pair of vectors $(\omega',\Omega)\in\mR^{n+N}$, $\beta'\in\mR^\nu$ is affinely $(\tau,\gamma,2)$-Diophantine.
\end{thm}

The set $\cT_\mu = \{\by=0, \, \bz=0\}$ is a reducible invariant $(n+N)$-torus of the system~\eqref{goal} and, consequently, of the system~\eqref{eqpert}. This torus is also invariant under the involution~\eqref{invol}. The torus $\cT_\mu$ is analytic and depends on $\mu\in\cG$ in a $C^\infty$-way in the sense of Whitney. It is characterized by frequency vector $(\omega',\Omega)\in\mR^{n+N}$ and Floquet matrix $\bfzero_m\oplus M'\in\gl(m+2p,\mR)$. The ball $\Gamma$ in Theorem~\ref{thmain} can be replaced by any smaller ball centered at the origin.

The informal meaning of Theorem~\ref{thmain} is obvious. If the perturbation terms $f$, $g$, $h$ are small enough then for most values of $\mu$ one can choose $\sigma$ in such a way that the system~\eqref{eqpert} will admit an invariant $(n+N)$-torus with frequency vector $(\omega',\Omega)$ and Floquet matrix $\bfzero_m\oplus M'$. Here $\omega'$ is close to $F(\mu)$ and $M'$ is close to $M(\mu)$.

\section{The Source Theorem in the Reversible Context~2}\label{source}

The dependence of the frequencies and nonzero Floquet exponents of the invariant torus $\{y=0, \, z=0\}$ of the unperturbed systems~\eqref{add3} on the parameter $\mu$ (for $\sigma=0$) is characterized by the mapping
\[
\mu \mapsto \bigl(F(\mu),\alpha(\mu),\beta(\mu)\bigr)\in\mR^{n+p}.
\]
According to the general ideology of Herman's method, the analogous mapping in the source theorem should be submersive. This is achieved by replacing $F(\mu)+\Delta(\sigma,\mu)$ with an independent external parameter $\omega\in\mR^n$ and assuming that the mapping $\mu \mapsto \bigl(\alpha(\omega,\mu),\beta(\omega,\mu)\bigr)\in\mR^p$ (for $M$ dependent on $\omega$) is submersive for fixed $\omega$.

To be more precise, let $n\in\mZ_+$, $m\in\mN$, $p\in\mZ_+$, $s\in\mZ_+$, and $\omega_\star\in\mR^n$. Consider an analytic $(m+n+s)$-parameter family of analytic differential equations
\begin{equation}
\begin{aligned}
\dot{x} &= \omega+\xi(y,z,\sigma,\omega,\mu)+f(x,y,z,\sigma,\omega,\mu), \\
\dot{y} &= \sigma+\eta(y,z,\sigma,\omega,\mu)+g(x,y,z,\sigma,\omega,\mu), \\
\dot{z} &= M(\omega,\mu)z+\zeta(y,z,\sigma,\omega,\mu)+h(x,y,z,\sigma,\omega,\mu),
\end{aligned}
\label{eqauto}
\end{equation}
where $x\in\mT^n$, $y\in\cO_m(0)$, $z\in\cO_{2p}(0)$ are the phase space variables, $\sigma\in\cO_m(0)$, $\omega\in\cO_n(\omega_\star)$, $\mu\in\cO_s(0)$ are external parameters, $M$ is a $2p\times 2p$ matrix-valued function, and $\xi=O(y,z)$, $\eta=O_2(y,z)$, $\zeta=O_2(y,z,\sigma)$. The functions $M$, $\xi$, $\eta$, $\zeta$ are supposed to be fixed whereas the summands $f$, $g$, $h$ are small perturbation terms. Let the systems~\eqref{eqauto} be reversible with respect to the phase space involution
\begin{equation}
G:(x,y,z)\mapsto(-x,-y,Rz)
\label{ution}
\end{equation}
(see~\eqref{rev2}), where $R\in\GL(2p,\mR)$ is an involutive matrix with eigenvalues $1$ and $-1$ of multiplicity $p$ each, $M(\omega,\mu)R\equiv-RM(\omega,\mu)$, and the spectrum of $M(\omega_\star,0)$ is simple. One may assume that the spectrum of $M(\omega,\mu)$ is simple for any $\omega$ and $\mu$ and has the form $\fS\bigl(\nu_1,\nu_2,\nu_3;\alpha(\omega,\mu),\beta(\omega,\mu)\bigr)$ where $\nu_1+\nu_2+2\nu_3=p$. Retain the notation $\nu=\nu_2+\nu_3\in\mZ_+$.

\begin{thm}[\cite{S16RCD}]\label{thsource}
Suppose that the mapping
\begin{equation}
\mu \mapsto \bigl(\alpha(\omega_\star,\mu),\beta(\omega_\star,\mu)\bigr)\in\mR^p
\label{mapmu}
\end{equation}
is \emph{submersive} at the origin $\mu=0$ \textup{(}so that $s\geq p$\textup{)}. Then there exists a neighborhood $\fO\subset\mR^{n+s}$ of the point $(\omega_\star,0)$ such that for any closed set $\Gamma\subset\fO$ that is diffeomorphic to an $(n+s)$-dimensional ball and contains the point $(\omega_\star,0)$ in its interior, the following holds. For every complex neighborhood
\[
\cC\subset(\mC/2\pi\mZ)^n\times\mC^{2m+2p+n+s}
\]
of the set
\[
\mT^n\times\{0\in\mR^m\}\times\{0\in\mR^{2p}\}\times\{0\in\mR^m\}\times\Gamma
\]
and every $\cL\in\mN$, $\vare>0$, $\tau>n-1$ \textup{(}$\tau\geq 0$ for $n=0$\textup{)}, $\gamma>0$, there is a number $\delta>0$ with the following properties.

Suppose that the perturbation terms $f$, $g$, $h$ in~\eqref{eqauto} can be holomorphically continued to the neighborhood $\cC$ and $|f|<\delta$, $|g|<\delta$, $|h|<\delta$ in $\cC$. Then for each $(\omega_0,\mu_0)\in\Gamma$, there exist points
\begin{equation}
v(\omega_0,\mu_0)\in\mR^m, \quad u(\omega_0,\mu_0)\in\mR^n, \quad w(\omega_0,\mu_0)\in\mR^s
\label{sdvig}
\end{equation}
and a change of variables
\begin{equation}
\begin{aligned}
x &= \bx+a(\bx,\omega_0,\mu_0), \\
y &= \by+b^0(\bx,\omega_0,\mu_0)+b^1(\bx,\omega_0,\mu_0)\by+b^2(\bx,\omega_0,\mu_0)\bz, \\
z &= \bz+c^0(\bx,\omega_0,\mu_0)+c^1(\bx,\omega_0,\mu_0)\by+c^2(\bx,\omega_0,\mu_0)\bz
\end{aligned}
\label{trans}
\end{equation}
\textup{(}cf.~\eqref{add4}\textup{)} with $\bx\in\mT^n$, $\by\in\cO_m(0)$, $\bz\in\cO_{2p}(0)$ such that the following holds.

\textup{i)} The coefficients $a$, $b^0$, $b^1$, $b^2$, $c^0$, $c^1$, $c^2$ in~\eqref{trans} are mappings ranging in $\mR^n$, $\mR^m$, $\gl(m,\mR)$, $\mR^{m\times 2p}$, $\mR^{2p}$, $\mR^{2p\times m}$, $\gl(2p,\mR)$, respectively. These mappings are analytic in $\bx$ and $C^\infty$-smooth in $(\omega_0,\mu_0)$. All the partial derivatives of each component of these mappings of any order from $0$ to $\cL$ are smaller than $\vare$ in absolute value everywhere in $\mT^n\times\Gamma$. The functions $u$, $v$, $w$ in~\eqref{sdvig} are $C^\infty$-smooth as functions in $(\omega_0,\mu_0)$, and all the partial derivatives of each component of these functions of any order from $0$ to $\cL$ are smaller than $\vare$ in absolute value everywhere in $\Gamma$. Moreover, $\{a\}_{\bx}\equiv 0$.

\textup{ii)} For each $(\omega_0,\mu_0)\in\Gamma$, the change of variables~\eqref{trans} commutes with the involution~\eqref{ution} in the sense that in the new variables $(\bx,\by,\bz)$, the involution $G$ takes the form
\[
G:(\bx,\by,\bz)\mapsto(-\bx,-\by,R\bz).
\]

\textup{iii)} For \emph{any} point $(\omega_0,\mu_0)\in\Gamma$ such that the pair of vectors $\omega_0\in\mR^n$, $\beta(\omega_0,\mu_0)\in\mR^\nu$ is affinely $(\tau,\gamma,2)$-Diophantine, the system~\eqref{eqauto} at the parameter values
\begin{equation}
\sigma=v(\omega_0,\mu_0), \quad \omega=\omega_0+u(\omega_0,\mu_0), \quad \mu=\mu_0+w(\omega_0,\mu_0)
\label{shift}
\end{equation}
after the coordinate transformation~\eqref{trans} takes the form
\begin{equation}
\dot{\bx}=\omega_0+O(\by,\bz), \quad \dot{\by}=O_2(\by,\bz), \quad \dot{\bz}=M(\omega_0,\mu_0)\bz+O_2(\by,\bz).
\label{ideal}
\end{equation}
\end{thm}

Consider any point $(\omega_0,\mu_0)\in\Gamma$ such that the pair of vectors $\omega_0\in\mR^n$, $\beta(\omega_0,\mu_0)\in\mR^\nu$ is affinely $(\tau,\gamma,2)$-Diophantine. The system~\eqref{eqauto} without the terms $f$, $g$, $h$ (the unperturbed system) at the parameter values $\sigma=0$, $\omega=\omega_0$, $\mu=\mu_0$ admits the reducible invariant $n$-torus $\{y=0, \, z=0\}$ with frequency vector $\omega_0\in\mR^n$ and Floquet matrix $\bfzero_m\oplus M(\omega_0,\mu_0)\in\gl(m+2p,\mR)$. According to the equations~\eqref{ideal}, the perturbed system~\eqref{eqauto} at the \emph{shifted} parameter values~\eqref{shift} has the reducible invariant $n$-torus $\{\by=0, \, \bz=0\}$ with \emph{the same} frequency vector and Floquet matrix. This torus is analytic and depends on $\omega_0$ and $\mu_0$ in a $C^\infty$-way in the sense of Whitney.

Theorem~\ref{thsource} is a particular case of the main result of \cite{S16RCD}, see the theorem in Section~2 of \cite{S16RCD} and the remark at the end of that section. Indeed, since the spectrum of the matrix $M(\omega_\star,0)$ is simple, the submersivity of the mapping~\eqref{mapmu} at $\mu=0$ is tantamount to the statement that $\mu\mapsto M(\omega_\star,\mu)$ is a versal unfolding of $M(\omega_\star,0)$ with respect to the adjoint action of the group of the $2p\times 2p$ matrices commuting with $R$ on the space of the $2p\times 2p$ matrices anti-commuting with $R$ \cite{H96JDE,S91TSP,S93CJM}. The paper \cite{S16RCD} treats the general case where $\det M(\omega_\star,0)\neq 0$ but $M(\omega_\star,0)$ is allowed to possess multiple eigenvalues (and even to be non-diagonalizable over $\mC$). The equality $\{a\}_{\bx}\equiv 0$ was not claimed in \cite{S16RCD} but the possibility of achieving it is obvious: if this equality is not valid, one should just replace $a$ with $a-\{a\}_{\bx}$. Besides, the main theorem of \cite{S16RCD} asserts the existence of a suitable set $\Gamma$ rather than the suitability of any sufficiently ``small'' $\Gamma$ containing the point $(\omega_\star,0)$ in its interior, but such a sharpening is also straightforward (see again the remark at the end of Section~2 of \cite{S16RCD}).

\section{A Proof of Theorem~\ref{thmain}}\label{proof}

Our goal is to deduce Theorem~\ref{thmain} from Theorem~\ref{thsource} following the general Herman-like scheme of \cite{S07DCDS}. Let the systems~\eqref{eqpert} satisfy the hypotheses of Theorem~\ref{thmain}. Since $M(\mu)$ depends on $\mu$ analytically and the spectrum of $M(0)$ is simple, one can introduce an additional parameter $\chi\in\cO_S(0)$ for an appropriate $S\in\mZ_+$ and construct an analytic family $M\new(\mu,\chi)$ of $2p\times 2p$ real matrices such that the following holds.

(a) $M\new(\mu,0)\equiv M(\mu)$ and $M\new(\mu,\chi)R\equiv-RM\new(\mu,\chi)$. As a consequence, one may assume that for any $\mu$ and $\chi$, the spectrum of $M\new(\mu,\chi)$ is simple and has the form
\[
\fS\bigl(\nu_1,\nu_2,\nu_3;\alpha\new(\mu,\chi),\beta\new(\mu,\chi)\bigr)
\]
where $\alpha\new(\mu,0)\equiv\alpha(\mu)$ and $\beta\new(\mu,0)\equiv\beta(\mu)$.

(b) The mapping
\[
(\mu,\chi) \mapsto \bigl(\alpha\new(\mu,\chi),\beta\new(\mu,\chi)\bigr)\in\mR^p
\]
is \emph{submersive} at $\mu=0$, $\chi=0$ (so that $s+S\geq p$).

The existence of a $2p\times 2p$ matrix-valued function $M\new$ satisfying these conditions follows immediately from the theory of normal forms and versal unfoldings of infinitesimally reversible matrices \cite{H96JDE,S91TSP,S93CJM}. It always suffices to set $S=p$.

Now introduce two more additional parameters $\omega\in\cO_n\bigl(F(0)\bigr)$ and $\theta\in\cO_N(0)$ and consider the analytic $(m+s+S+n+N)$-parameter family of analytic differential equations
\begin{equation}
\begin{aligned}
\dot{x} &= \omega+\xi(y,z,\sigma,\mu)+f(x,y,z,\sigma,\mu,X), \\
\dot{y} &= \sigma+\eta(y,z,\sigma,\mu)+g(x,y,z,\sigma,\mu,X), \\
\dot{z} &= M\new(\mu,\chi)z+\zeta(y,z,\sigma,\mu)+h(x,y,z,\sigma,\mu,X), \\
\dot{X} &= \Omega+\theta.
\end{aligned}
\label{eqext}
\end{equation}
The systems~\eqref{eqext} are reversible with respect to the involution~\eqref{invol} and satisfy all the hypotheses of Theorem~\ref{thsource}, with

(1) $n+N$, $(x,X)$, $(\omega,\Omega+\theta)$, $\bigl(F(0),\Omega\bigr)$ playing the roles of $n$, $x$, $\omega$, $\omega_\star$, respectively;

(2) $s+S$, $(\mu,\chi)$, $M\new$ playing the roles of $s$, $\mu$, $M$, respectively;

(3) $\fG$ playing the role of $G$.

Theorem~\ref{thsource} provides us with closed balls $\Gamma\subset\mR^s$, $\Gamma_1\subset\mR^S$, $\Gamma_2\subset\mR^n$, $\Gamma_3\subset\mR^N$ centered at $0$, $0$, $F(0)$, $0$, respectively (the radii of these balls can be chosen to be arbitrarily small), such that for every complex neighborhood~\eqref{neigh} of the set~\eqref{set} and every $\cL\in\mN$, $\tau>n+N-1$, $\gamma>0$, the following holds. Suppose that the perturbation terms $f$, $g$, $h$ in~\eqref{eqpert} (and~\eqref{eqext}) can be holomorphically continued to the neighborhood~\eqref{neigh} and are sufficiently small in~\eqref{neigh}. Then for any $\mu_0\in\Gamma$, $\chi_0\in\Gamma_1$, $\omega_0\in\Gamma_2$, $\theta_0\in\Gamma_3$, there exist points
\begin{equation}
\begin{gathered}
v(\omega_0,\theta_0,\mu_0,\chi_0)\in\mR^m, \quad u(\omega_0,\theta_0,\mu_0,\chi_0)\in\mR^n, \quad U(\omega_0,\theta_0,\mu_0,\chi_0)\in\mR^N, \\
w(\omega_0,\theta_0,\mu_0,\chi_0)\in\mR^s, \quad W(\omega_0,\theta_0,\mu_0,\chi_0)\in\mR^S
\end{gathered}
\label{bigsdvig}
\end{equation}
and a change of variables
\begin{equation}
\begin{aligned}
x &= \bx+a(\bx,\bX,\omega_0,\theta_0,\mu_0,\chi_0), \\
y &= \by+b^0(\bx,\bX,\omega_0,\theta_0,\mu_0,\chi_0)+b^1(\bx,\bX,\omega_0,\theta_0,\mu_0,\chi_0)\by+b^2(\bx,\bX,\omega_0,\theta_0,\mu_0,\chi_0)\bz, \\
z &= \bz+c^0(\bx,\bX,\omega_0,\theta_0,\mu_0,\chi_0)+c^1(\bx,\bX,\omega_0,\theta_0,\mu_0,\chi_0)\by+c^2(\bx,\bX,\omega_0,\theta_0,\mu_0,\chi_0)\bz, \\
X &= \bX+A(\bx,\bX,\omega_0,\theta_0,\mu_0,\chi_0)
\end{aligned}
\label{bigtrans}
\end{equation}
with $\bx\in\mT^n$, $\bX\in\mT^N$, $\by\in\cO_m(0)$, $\bz\in\cO_{2p}(0)$ such that the following is valid.

First, the coefficients $a$, $A$, $b^0$, $b^1$, $b^2$, $c^0$, $c^1$, $c^2$ in~\eqref{bigtrans} are analytic in $(\bx,\bX)$ and $C^\infty$-smooth in $(\omega_0,\theta_0,\mu_0,\chi_0)$. The functions $u$, $U$, $v$, $w$, $W$ in~\eqref{bigsdvig} are $C^\infty$-smooth. All the mappings $a$, $A$, $b^0$, $b^1$, $b^2$, $c^0$, $c^1$, $c^2$, $u$, $U$, $v$, $w$, $W$ are small in the $C^{\cL}$-topology. Moreover, $\{a\}_{(\bx,\bX)}\equiv 0$ and $\{A\}_{(\bx,\bX)}\equiv 0$.

Second, for any $\mu_0\in\Gamma$, $\chi_0\in\Gamma_1$, $\omega_0\in\Gamma_2$, $\theta_0\in\Gamma_3$, the change of variables~\eqref{bigtrans} commutes with the involution~\eqref{invol}.

Third, for \emph{any} points $\mu_0\in\Gamma$, $\chi_0\in\Gamma_1$, $\omega_0\in\Gamma_2$, $\theta_0\in\Gamma_3$ such that the pair of vectors
\begin{equation}
(\omega_0,\Omega+\theta_0)\in\mR^{n+N}, \quad \beta\new(\mu_0,\chi_0)\in\mR^\nu
\label{vectors2}
\end{equation}
is affinely $(\tau,\gamma,2)$-Diophantine, the system~\eqref{eqext} at the parameter values
\begin{equation}
\begin{gathered}
\sigma=v(\omega_0,\theta_0,\mu_0,\chi_0), \quad \omega=\omega_0+u(\omega_0,\theta_0,\mu_0,\chi_0), \quad \theta=\theta_0+U(\omega_0,\theta_0,\mu_0,\chi_0), \\
\mu=\mu_0+w(\omega_0,\theta_0,\mu_0,\chi_0), \quad \chi=\chi_0+W(\omega_0,\theta_0,\mu_0,\chi_0)
\end{gathered}
\label{bigshift}
\end{equation}
after the coordinate transformation~\eqref{bigtrans} takes the form
\begin{equation}
\begin{gathered}
\dot{\bx}=\omega_0+O(\by,\bz), \quad \dot{\by}=O_2(\by,\bz), \\
\dot{\bz}=M\new(\mu_0,\chi_0)\bz+O_2(\by,\bz), \quad \dot{\bX}=\Omega+\theta_0+O(\by,\bz).
\end{gathered}
\label{bigideal}
\end{equation}

We claim that $U(\omega_0,\theta_0,\mu_0,\chi_0)=0$ and $A(\bx,\bX,\omega_0,\theta_0,\mu_0,\chi_0)\equiv 0$ whenever the pair of vectors~\eqref{vectors2} is affinely $(\tau,\gamma,2)$-Diophantine (the shifts along $X$ and $\theta$ vanish due to the special form of the equation for $\dot{X}$ in~\eqref{eqext}). Indeed, at the parameter values~\eqref{bigshift} one has
\[
\dot{X} = \Omega+\theta_0+U
\]
according to~\eqref{eqext}, where $U=U(\omega_0,\theta_0,\mu_0,\chi_0)$. On the other hand,
\[
\dot{X} = \Omega+\theta_0+O(\by,\bz)
+ \frac{\partial A}{\partial\bx}\bigl[\omega_0+O(\by,\bz)\bigr]
+ \frac{\partial A}{\partial\bX}\bigl[\Omega+\theta_0+O(\by,\bz)\bigr]
\]
in virtue of~\eqref{bigtrans} and~\eqref{bigideal}, where $A=A(\bx,\bX,\omega_0,\theta_0,\mu_0,\chi_0)$. For $\by=0$, $\bz=0$ we get
\[
U = \frac{\partial A}{\partial\bx}\omega_0+\frac{\partial A}{\partial\bX}(\Omega+\theta_0).
\]
Since the components of the vector $(\omega_0,\Omega+\theta_0)$ are rationally independent and $\{A\}_{(\bx,\bX)}=0$, we arrive at the claim (cf.\ \cite[section~6.2]{S07DCDS}). Thus, one can set
\[
U(\omega,\theta,\mu,\chi)\equiv 0, \quad A(\bx,\bX,\omega,\theta,\mu,\chi)\equiv 0,
\]
and $\bX=X$ in~\eqref{bigtrans}. In the sequel, we will also always set $\theta_0=0$.

One may assume that $F(\Gamma)$ lies in the interior of $\Gamma_2$. If the functions $u$, $v$, $w$, $W$ are small enough, then the system of equations
\begin{equation}
\begin{aligned}
\omega+u(\omega,0,\mu,\chi) &= F\bigl(\mu+w(\omega,0,\mu,\chi)\bigr)+\Delta\bigl(v(\omega,0,\mu,\chi), \, \mu+w(\omega,0,\mu,\chi)\bigr), \\
\chi+W(\omega,0,\mu,\chi) &= 0
\end{aligned}
\label{keysystem}
\end{equation}
with $\mu\in\Gamma$ can be solved with respect to $\omega$ and $\chi$:
\[
\omega=\varp(\mu), \quad \chi=\psi(\mu),
\]
where $\varp:\Gamma\to\Gamma_2$ and $\psi:\Gamma\to\Gamma_1$ are $C^\infty$-functions close to $F$ and $0$, respectively, in the $C^{\cL}$-topology. The key observation is that for any $\mu_0\in\Gamma$, the system~\eqref{eqext} at the parameter values~\eqref{bigshift} with $\omega_0=\varp(\mu_0)$, $\chi_0=\psi(\mu_0)$, $\theta_0=0$ coincides with the \emph{original} system~\eqref{eqpert} at the parameter values
\begin{equation}
\sigma=v(\omega_0,0,\mu_0,\chi_0), \quad \mu=\mu_0+w(\omega_0,0,\mu_0,\chi_0)
\label{sigmamu}
\end{equation}
(recall that $U\equiv 0$). Indeed, if $\omega_0=\varp(\mu_0)$, $\chi_0=\psi(\mu_0)$, $\theta_0=0$ then the equations~\eqref{keysystem} imply that the values of the parameters $\sigma$, $\omega$, $\mu$, $\chi$, $\theta$ given by~\eqref{bigshift} satisfy the relations
\[
\omega=F(\mu)+\Delta(\sigma,\mu), \quad \chi=0, \quad \theta=0.
\]

Let $\vare_2\in(0,1)$ be an arbitrary number, and let $\Gamma'\subset\Gamma$ be the closed ball centered at the origin and determined by the equality
\begin{equation}
\meas_s(\Gamma\setminus\Gamma') = \frac{\vare_2}{3}\meas_s\Gamma.
\label{meas1}
\end{equation}
If the functions $u$, $v$, $w$, $W$ are small enough, then the equation
\begin{equation}
\mu=\mu_0+w\bigl(\varp(\mu_0),0,\mu_0,\psi(\mu_0)\bigr)
\label{keyeq}
\end{equation}
with $\mu\in\Gamma'$ can be solved with respect to $\mu_0$:
\[
\mu_0=\Upsilon(\mu),
\]
where $\Upsilon:\Gamma'\to\Gamma$ is a $C^\infty$-function close to the identity mapping $\mu\mapsto\mu$ in the $C^{\cL}$-topology.

Introduce the functions
\begin{align*}
\Phi(\mu) &= \varp\bigl(\Upsilon(\mu)\bigr), \\
\Psi(\mu) &= \psi\bigl(\Upsilon(\mu)\bigr), \\
\Theta(\mu) &= v\bigl(\Phi(\mu),0,\Upsilon(\mu),\Psi(\mu)\bigr)
\end{align*}
for $\mu\in\Gamma'$. The functions $\Phi:\Gamma'\to\Gamma_2$, $\Psi:\Gamma'\to\Gamma_1$, $\Theta:\Gamma'\to\mR^m$ are well-defined, $C^\infty$-smooth, and close to $F$, $0$, $0$, respectively, in the $C^{\cL}$-topology provided that $u$, $v$, $w$, $W$ are small enough. Note that $\Theta$ can be continued to a small (in the $C^{\cL}$-topology) $C^\infty$-function $\Theta:\Gamma\to\mR^m$ (see \cite[\S~6.1.4]{BHS96LNM} and references therein).

We arrive at the following conclusion. For any point $\mu\in\Gamma'$, set
\[
\mu_0=\Upsilon(\mu), \quad \omega_0=\Phi(\mu), \quad \chi_0=\Psi(\mu), \quad \theta_0=0.
\]
If the pair of vectors
\begin{equation}
(\omega_0,\Omega)\in\mR^{n+N}, \quad \beta\new(\mu_0,\chi_0)\in\mR^\nu
\label{vectors3}
\end{equation}
is affinely $(\tau,\gamma,2)$-Diophantine, then the \emph{original} system~\eqref{eqpert} at the parameter values $\mu$ and $\sigma=\Theta(\mu)$ (compare~\eqref{sigmamu} and~\eqref{keyeq}) after the $\fG$-commuting coordinate transformation~\eqref{bigtrans} with $A\equiv 0$ and $\bX=X$ takes the form~\eqref{goal} with
\[
\omega'=\omega_0, \quad M'=M\new(\mu_0,\chi_0).
\]
As in the case of the function $\Theta$, one may regard the transformation~\eqref{change} defined this way as dependent on $\mu\in\Gamma$ rather than on $\mu\in\Gamma'$.

It remains to estimate the measure of the set of points $\mu\in\Gamma'$ for which the pair of vectors~\eqref{vectors3} is affinely $(\tau,\gamma,2)$-Diophantine. Suppose that $\tau>(n+N)Q$, $\tau\geq\tau_\ast$, and $\cL\geq Q$. Let $\Gamma''\subset\Gamma'$ be the closed ball centered at the origin and determined by the equality
\begin{equation}
\meas_s(\Gamma'\setminus\Gamma'') = \frac{\vare_2}{3}\meas_s\Gamma.
\label{meas2}
\end{equation}
If the radius of $\Gamma$ is equal to $r$ then the relations~\eqref{meas1} and~\eqref{meas2} mean that the radii of $\Gamma'$ and $\Gamma''$ are equal to $r(1-\vare_2/3)^{1/s}$ and $r(1-2\vare_2/3)^{1/s}$, respectively. One may assume that the pair of mappings $F$, $\beta$ is affinely $(Q,2)$-nondegenerate at each point of $\Gamma''$. Apply Lemma~\ref{lmDioph} to the case where $\fK$ is the interior of $\Gamma'$, $K=\Gamma''$ and
\[
\tF(\mu)=\Phi(\mu), \quad \tbeta(\mu)=\beta\new\bigl(\Upsilon(\mu),\Psi(\mu)\bigr).
\]
The functions $\tF$ and $\tbeta$ are $C^\infty$-smooth and $C^{\cL}$-close in $\Gamma'$ to $F$ and $\beta$, respectively.

According to Lemma~\ref{lmDioph}, there exists a number $\gamma_0=\gamma_0(\vare_2,\tau,\gamma_\ast)>0$ such that the following holds. Let $\tF$ and $\tbeta$ be sufficiently close in $\Gamma'$ to $F$ and $\beta$, respectively, in the $C^Q$-topology. Then for any $\gamma\in(0,\gamma_0]$ the measure of the set $\cG$ of those points $\mu\in\Gamma''$ for which the pair of vectors~\eqref{vectors1} is affinely $(\tau,\gamma,2)$-Diophantine, is greater than
\[
\frac{3(1-\vare_2)}{3-2\vare_2}\meas_s\Gamma''.
\]
Thus,
\begin{equation}
\meas_s(\Gamma''\setminus\cG) < \frac{\vare_2}{3-2\vare_2}\meas_s\Gamma'' = \frac{\vare_2}{3}\meas_s\Gamma.
\label{meas3}
\end{equation}
Combining~\eqref{meas1}, \eqref{meas2}, and~\eqref{meas3}, one obtains $\meas_s(\Gamma\setminus\cG) < \vare_2\meas_s\Gamma$. This completes the proof of Theorem~\ref{thmain}.


\begin{thebibliography}{99}

\bibitem{AL12JDDE}
H.~N.~Alishah and R.~de la Llave, \textit{Tracing KAM tori in presymplectic dynamical systems}, J. Dynam. Differential Equations \textbf{24} (2012), no.~4, 685--711. MR 3000600

\bibitem{A84}
V.~I.~Arnold, \textit{Reversible systems}, Nonlinear and turbulent processes in physics, Vol.~3 (Kiev, 1983), Harwood Academic Publ., Chur, 1984, pp.~1161--1174. MR 0824779

\bibitem{AKN06}
V.~I.~Arnold, V.~V.~Kozlov, and A.~I.~Neishtadt, \textit{Mathematical aspects of classical and celestial mechanics}, 3rd ed., Encyclopaedia of Mathematical Sciences, vol.~3, Springer-Verlag, Berlin, 2006. MR 2269239

\bibitem{AS86}
V.~I.~Arnold and M.~B.~Sevryuk, \textit{Oscillations and bifurcations in reversible systems}, Nonlinear phenomena in plasma physics and hydrodynamics (R.~Z.~Sagdeev, ed.), Mir, Moscow, 1986, pp.~31--64.

\bibitem{BCHV09PD}
H.~W.~Broer, M.~C.~Ciocci, H.~Han{\ss}mann, and A.~Vanderbauwhede, \textit{Quasi-periodic stability of normally resonant tori}, Phys.~D \textbf{238} (2009), no.~3, 309--318. MR 2590451

\bibitem{BHN07JDE}
H.~W.~Broer, J.~Hoo, and V.~Naudot, \textit{Normal linear stability of quasi-periodic tori}, J. Differential Equations \textbf{232} (2007), no.~2, 355--418. MR 2286385

\bibitem{BH95JDDE}
H.~W.~Broer and G.~B.~Huitema, \textit{Unfoldings of quasi-periodic tori in reversible systems}, J. Dynam. Differential Equations \textbf{7} (1995), no.~1, 191--212. MR 1321710

\bibitem{BHS96G}
H.~W.~Broer, G.~B.~Huitema, and M.~B.~Sevryuk, \textit{Families of quasi-periodic motions in dynamical systems depending on parameters}, Nonlinear dynamical systems and chaos (Groningen, 1995), Progr. Nonlinear Differential Equations Appl., vol.~19, Birkh\"auser, Basel, 1996, pp.~171--211. MR 1391497

\bibitem{BHS96LNM}
H.~W.~Broer, G.~B.~Huitema, and M.~B.~Sevryuk, \textit{Quasi-periodic motions in families of dynamical systems. Order amidst chaos}, Lecture Notes in Mathematics, vol.~1645, Springer-Verlag, Berlin, 1996. MR 1484969

\bibitem{BHTB90}
H.~W.~Broer, G.~B.~Huitema, F.~Takens, and B.~L.~J.~Braaksma, \textit{Unfoldings and bifurcations of quasi-periodic tori}, Mem. Amer. Math. Soc. \textbf{83} (1990), no.~421. MR 1041003

\bibitem{BS10}
H.~W.~Broer and M.~B.~Sevryuk, \textit{KAM theory:\ Quasi-periodicity in dynamical systems}, Handbook of Dynamical Systems, Vol.~3 (H.~W.~Broer, B.~Hasselblatt, and F.~Takens, eds.), Elsevier~B.V., Amsterdam, 2010, pp.~249--344.

\bibitem{CCL13JDE}
R.~C.~Calleja, A.~Celletti, and R.~de la Llave, \textit{A KAM theory for conformally symplectic systems:\ efficient algorithms and their validation}, J. Differential Equations \textbf{255} (2013), no.~5, 978--1049. MR 3062760

\bibitem{D76TAMS}
R.~L.~Devaney, \textit{Reversible diffeomorphisms and flows}, Trans. Amer. Math. Soc. \textbf{218} (1976), 89--113. MR 0402815

\bibitem{FW16}
A.~Fortunati and S.~Wiggins, \textit{The Kolmogorov--Arnold--Moser \textup{(}KAM\textup{)} and Nekhoroshev theorems with arbitrary time dependence}, Essays in mathematics and its applications. In honor of Vladimir Arnold (T.~M.~Rassias and P.~M.~Pardalos, eds.), Springer-Verlag, Cham, 2016, pp.~89--99. MR 3526916

\bibitem{GHL14}
A.~Gonz\'alez-Enr\'{\i}quez, \`A.~Haro, and R.~de la Llave, \textit{Singularity theory for non-twist KAM tori}, Mem. Amer. Math. Soc. \textbf{227} (2014), no.~1067. MR 3154587

\bibitem{HCFLM16}
\`A.~Haro, M.~Canadell, J.-L.~Figueras, A.~Luque, and J.-M.~Mondelo, \textit{The parameterization method for invariant manifolds. From rigorous results to effective computations}, Applied Mathematical Sciences, vol.~195, Springer-Verlag, Cham, 2016. MR 3467671

\bibitem{H96JDE}
I.~Hoveijn, \textit{Versal deformations and normal forms for reversible and Hamiltonian linear systems}, J. Differential Equations \textbf{126} (1996), no.~2, 408--442. MR 1383984

\bibitem{KMS16JMPA}
A.~Kiesenhofer, E.~Miranda, and G.~Scott, \textit{Action-angle variables and a KAM theorem for $b$-Poisson manifolds}, J. Math. Pures Appl. (9) \textbf{105} (2016), no.~1, 66--85. MR 3427939

\bibitem{KX14AMC}
Y.~Kong and J.~Xu, \textit{Persistence of lower dimensional hyperbolic tori for reversible system}, Appl. Math. Comput. \textbf{236} (2014), 408--421. MR 3197738

\bibitem{KP01UMZ}
A.~A.~Kubichka and I.~O.~Parasyuk, \textit{Bifurcation of a Whitney smooth family of coisotropic invariant tori of a Hamiltonian system under a small deformation of the symplectic structure}, Ukra\"{\i}n. Mat. Zh. \textbf{53} (2001), no.~5, 610--624 (Ukrainian). MR 1854551.
English translation: Ukrainian Math. J. \textbf{53} (2001), no.~5, 701--718.

\bibitem{LR98PD}
J.~S.~W.~Lamb and J.~A.~G.~Roberts, \textit{Time-reversal symmetry in dynamical systems:\ a survey}, Phys.~D \textbf{112} (1998), no.~1--2, 1--39. MR 1605826

\bibitem{L71DAN}
V.~F.~Lazutkin, \textit{On the asymptotics of the eigenfunctions of the Laplace operator}, Dokl. Akad. Nauk SSSR \textbf{200} (1971), no.~6, 1277--1279 (Russian). MR 0310466.
English translation: Sov. Math., Dokl. \textbf{12} (1971), 1569--1572.

\bibitem{LY02ETDS}
Y.~Li and Y.~Yi, \textit{Persistence of invariant tori in generalized Hamiltonian systems}, Ergodic Theory Dynam. Systems \textbf{22} (2002), no.~4, 1233--1261. MR 1926285

\bibitem{LZH06JMAA}
B.~Liu, W.~Zhu, and Y.~Han, \textit{Persistence of lower-dimensional hyperbolic invariant tori for generalized Hamiltonian systems}, J. Math. Anal. Appl. \textbf{322} (2006), no.~1, 251--275. MR 2239236

\bibitem{LP05NK}
Yu.~V.~Love\u{\i}kin and I.~O.~Parasyuk, \textit{A theorem on a perturbation of coisotropic invariant tori of locally Hamiltonian systems and its applications}, Nelini\u{\i}ni Kolyv. \textbf{8} (2005), no.~4, 490--515 (Ukrainian). MR 2230790.
English translation: Nonlinear Oscil. (N.\,Y.) \textbf{8} (2005), no.~4, 487--512.

\bibitem{LP07UMZ}
Yu.~V.~Love\u{\i}kin and I.~O.~Parasyuk, \textit{Invariant tori of locally Hamiltonian systems that are close to conditionally integrable systems}, Ukra\"{\i}n. Mat. Zh. \textbf{59} (2007), no.~1, 71--98 (Ukrainian). MR 2356031.
English translation: Ukrainian Math. J. \textbf{59} (2007), no.~1, 70--99.

\bibitem{M67MA}
J.~Moser, \textit{Convergent series expansions for quasi-periodic motions}, Math. Ann. \textbf{169} (1967), no.~1, 136--176. MR 0208078

\bibitem{P01}
J.~P\"oschel, \textit{A lecture on the classical KAM theorem}, Smooth ergodic theory and its applications (Seattle, WA, 1999), Proc. Sympos. Pure Math., vol.~69, Amer. Math. Soc., Providence, RI, 2001, pp.~707--732. MR 1858551

\bibitem{RQ92PR}
J.~A.~G.~Roberts and G.~R.~W.~Quispel, \textit{Chaos and time-reversal symmetry. Order and chaos in reversible dynamical systems}, Phys. Rep. \textbf{216} (1992), no.~2--3, 63--177. MR 1173588

\bibitem{S86}
M.~B.~Sevryuk, \textit{Reversible systems}, Lecture Notes in Mathematics, vol.~1211, Springer-Verlag, Berlin, 1986. MR 0871875

\bibitem{S91TSP}
M.~B.~Sevryuk, \textit{Reversible linear systems and their versal deformations}, Trudy Sem. Petrovsk. No.~15 (1991), 33--54, 235 (Russian). MR 1294389.
English translation: J. Soviet Math. \textbf{60} (1992), no.~5, 1663--1680. MR 1181098

\bibitem{S94}
M.~B.~Sevryuk, \textit{New results in the reversible KAM theory}, Seminar on dynamical systems (St. Petersburg, 1991), Progr. Nonlinear Differential Equations Appl., vol.~12, Birkh\"auser, Basel, 1994, pp.~184--199. MR 1279398

\bibitem{S95C}
M.~B.~Sevryuk, \textit{The iteration-approximation decoupling in the reversible KAM theory}, Chaos \textbf{5} (1995), no.~3, 552--565. MR 1350654

\bibitem{S97}
M.~B.~Sevryuk, \textit{Excitation of elliptic normal modes of invariant tori in Hamiltonian systems}, Topics in singularity theory, Amer. Math. Soc. Transl. Ser.~2, vol.~180, Amer. Math. Soc., Providence, RI, 1997, pp.~209--218. MR 1767126

\bibitem{S01}
M.~B.~Sevryuk, \textit{Excitation of elliptic normal modes of invariant tori in volume preserving flows}, Global analysis of dynamical systems, Inst. Phys., Bristol, 2001, pp.~339--352. MR 1858482

\bibitem{S06N}
M.~B.~Sevryuk, \textit{Partial preservation of frequencies in KAM theory}, Nonlinearity \textbf{19} (2006), no.~5, 1099--1140. MR 2221801

\bibitem{S07DCDS}
M.~B.~Sevryuk, \textit{Invariant tori in quasi-periodic non-autonomous dynamical systems via Herman's method}, Discrete Contin. Dyn. Syst. \textbf{18} (2007), no.~2--3, 569--595. MR 2291912

\bibitem{S07TMIS}
M.~B.~Sevryuk, \textit{Partial preservation of frequencies and Floquet exponents in KAM theory}, Tr. Mat. Inst. Steklova \textbf{259} (2007), 174--202 (Russian). MR 2433684.
English translation: Proc. Steklov Inst. Math. \textbf{259} (2007), 167--195.

\bibitem{S08N}
M.~B.~Sevryuk, \textit{KAM tori:\ persistence and smoothness}, Nonlinearity \textbf{21} (2008), no.~10, T177--T185. MR 2439472

\bibitem{S11RCD}
M.~B.~Sevryuk, \textit{The reversible context~2 in KAM theory:\ the first steps}, Regul. Chaotic Dyn. \textbf{16} (2011), no.~1--2, 24--38. MR 2774376

\bibitem{S12MMJ}
M.~B.~Sevryuk, \textit{KAM theory for lower dimensional tori within the reversible context~2}, Mosc. Math. J. \textbf{12} (2012), no.~2, 435--455, 462. MR 2978764

\bibitem{S12IM}
M.~B.~Sevryuk, \textit{Quasi-periodic perturbations within the reversible context~2 in KAM theory}, Indag. Math. (N.\,S.) \textbf{23} (2012), no.~3, 137--150. MR 2948617

\bibitem{S16RCD}
M.~B.~Sevryuk, \textit{Whitney smooth families of invariant tori within the reversible context~2 of KAM theory}, Regul. Chaotic Dyn. \textbf{21} (2016), no.~6, 599--620. MR 3583939

\bibitem{S93CJM}
C.~W.~Shih, \textit{Normal forms and versal deformations of linear involutive dynamical systems}, Chinese J. Math. \textbf{21} (1993), no.~4, 333--347. MR 1247555

\bibitem{W10DCDSS}
F.~Wagener, \textit{A parametrised version of Moser's modifying terms theorem}, Discrete Contin. Dyn. Syst. Ser.~S \textbf{3} (2010), no.~4, 719--768. MR 2684071

\bibitem{WXZ10DCDSB}
X.~Wang, J.~Xu, and D.~Zhang, \textit{Persistence of lower dimensional elliptic invariant tori for a class of nearly integrable reversible systems}, Discrete Contin. Dyn. Syst. Ser.~B \textbf{14} (2010), no.~3, 1237--1249. MR 2670193

\bibitem{WXZ17DCDS}
X.~Wang, J.~Xu, and D.~Zhang, \textit{A KAM theorem for the elliptic lower dimensional tori with one normal frequency in reversible systems}, Discrete Contin. Dyn. Syst. \textbf{37} (2017), no.~4, 2141--2160. MR 3640592

\bibitem{Y92A}
J.-C.~Yoccoz, \textit{Travaux de Herman sur les tores invariants}, S\'eminaire Bourbaki, Vol.~1991/92, Ast\'erisque No.~206 (1992), Exp. No.~754, 4, 311--344. MR 1206072

\bibitem{ZC10FPTA}
D.~Zhang and R.~Cheng, \textit{On invariant tori of nearly integrable Hamiltonian systems with quasiperiodic perturbation}, Fixed Point Theory Appl. 2010, Art. ID 697343, 17~pp. MR 2735724

\end{thebibliography}
\end{document}